\documentclass[12pt]{amsart}

\usepackage{tikz-cd} 

\usepackage[english]{babel}
\usepackage[T1]{fontenc}
\usepackage{amsmath}
\usepackage{amstext}
\usepackage{amsfonts}
\usepackage{amssymb}
\usepackage{amsthm}
\usepackage{mathtools}
\usepackage{dsfont}
\usepackage{color}
\usepackage{graphicx}
\usepackage[colorinlistoftodos]{todonotes}
\usepackage[colorlinks=true, allcolors=blue]{hyperref}
\usepackage{float}
\usepackage[colorinlistoftodos]{todonotes}
\usepackage{theoremref}
\usepackage{enumitem}
\usepackage{dirtytalk}
\usepackage{todonotes}
\usepackage{marginnote} 

\usepackage{multicol} 
\usepackage{titletoc}

\allowdisplaybreaks
\let\oldaddcontentsline\addcontentsline

\newcommand{\starttocentries}{\let\addcontentsline\oldaddcontentsline}
\newtheorem{theorem}{Theorem}[section]
\newtheorem{lemma}[theorem]{Lemma}
\newtheorem{prop}[theorem]{Proposition}
\newtheorem{cor}[theorem]{Corollary}

\newtheorem*{cor*}{Corollary}
\newtheorem*{conjecture*}{Conjecture}
\newtheorem*{thm*}{Theorem}
\newtheorem*{lem*}{Lemma}

\newtheorem*{prop*}{Proposition}
\theoremstyle{definition}

\newtheorem*{defn*}{Definition}

\theoremstyle{remark}
\newtheorem{remark}[theorem]{Remark}

\newcommand{\cB}{\mathcal{B}}

\newcommand{\cD}{\mathcal{D}}
\newcommand{\cE}{\mathcal{E}}
\newcommand{\cF}{\mathcal{F}}

\newcommand{\cO}{\mathcal{O}}

\newcommand{\cQ}{\mathcal{Q}}

\newcommand{\cU}{\mathcal{U}}

\newcommand{\id}{\operatorname{id}}

\makeatletter \def\l@subsection{\@tocline{2}{0pt}{1pc}{5pc}{}} \def\l@subsection{\@tocline{2}{0pt}{2pc}{6pc}{}} \makeatother

\title{
On the inclusion $\cO_2 \subset \cQ_2$}
\author{Jacopo Bassi, Roberto Conti }

\begin{document}

\maketitle

\begin{abstract}
The diadic $C^*$-algebra $\cQ_2$ contains canonically a copy of the Cuntz algebra $\cO_2$. It is shown that the inclusion $\cO_2 \subset \cQ_2$ is $C^*$-irreducible and rigid. It follows that the injective envelopes of these two $C^*$-algebras are $*$-isomorphic.
\end{abstract}


\section{Introduction}

Following an influential paper by R{\o}rdam \cite{Ro}, in recent times there has been a widespread interest in detecting so-called $C^*$-irreducible inclusions of $C^*$-algebras. Recall that a (unital) inclusion $A \subset B$ of $C^*$-algebras is said to be 
 $C^*$-{\it irreducible} if all intermediate $C^*$-algebras $C$ sitting between $A$ and $B$ are simple. In particular, $A$ and $B$ themselves must be simple as well. Several examples of such inclusions have been discussed in the literature, see e.g. \cite{AmKa, AmUr, BeOm,LiSc,HR}.

The diadic $C^*$-algebra $\cQ_2$ has been investigated in detail in \cite{LarsenLi,ACR}. It is simple, nuclear and purely infinite with $K_0(\cQ_2)$ and $K_1(\cQ_2)$ both isomorphic to $\mathbb Z$. Moreover, it contains canonically a copy of the Cuntz algebra $\cO_2$ (introduced in \cite{Cu}), which is also simple, nuclear and purely infinite, but with trivial $K$-theory. 
One of the main goals 
of this short note is to prove that the obtained inclusion $\cO_2 \subset \cQ_2$ is $C^*$-irreducible, a fact that seems to have been unnoticed so far. 
We remark that we have not been able to exhibit/identify any nontrivial $C^*$-subalgebra sitting between $\mathcal{O}_2$ and $\mathcal{Q}_2$. We will comment a bit more on this in Section 2. 
As a matter of fact 
one might even wonder whether
the inclusion $\cO_2 \subset \cQ_2$ is {\it tight}, meaning that there is no nontrivial intermediate $C^*$-algebra, 
but 
for the time being 
we will not discuss this issue any further.
Nevertheless we also show that $\cO_2 \subset \cQ_2$ is {\it rigid}, namely the only ucp map of $\cQ_2$ that restricts to the identity of $\cO_2$ is the identity of $\cQ_2$. 
As a consequence, the injective envelopes of $\cO_2$ and $\cQ_2$ are isomorphic.

\section{Main results
}
The Cuntz algebra $\cO_2$ is the universal $C^*$-algebra generated by two isometries $S_1$ and $S_2$ satisfying 
\begin{equation*}S_1 S_1^* + S_2 S_2^* =1.
\end{equation*}
The diadic $C^*$-algebra $\cQ_2$ is the universal $C^*$-algebra generated by a unitary $U$ and an isometry $S_2$ such that 
\[U^2 S_2 = S_2 U , \quad S_2 S_2^* + U S_2 S_2^* U^* = 1\]
Setting $S_1 := U S_2$, it readily follows that $S_1$ and $S_2$ satisfy the relations of the generating isometries of $\cO_2$. This means that there is a canonical injection of $\cO_2$ inside $\cQ_2$.
As stated above, it is well-known that both $\cO_2$ and $\cQ_2$ are simple $C^*$-algebras, and moreover the so-called diagonal subalgebra $\cD_2 \subset \cO_2$ is Cartan not only in $\cO_2$ but also in $\cQ_2$.

Recall that given an inclusion of unital $C^*$-algebras $A \subset B$, a pseudo-expectation from $B$ to $A$ is a ucp map $\psi : B \rightarrow I(A)$ such that $\psi|_A = \id_A$, where $I(A)$ is the injective envelope of $A$ (\cite{ham} Definition 2.2); moreover, the inclusion is hereditarily essential if for every intermediate unital $C^*$-algebra $A \subset C \subset B$, every nonzero ideal in $C$ intersects $A$ non-trivially (\cite{PiZa} Definition 3.4).

\begin{lemma}
\label{lem1}
    Let $A \subset B$ be a unital inclusion of unital simple $C^*$-algebras. Then $A \subset B$ is hereditarily essential if and only if it is $C^*$-irreducible.
\end{lemma}
\begin{proof}
    Suppose that $A \subset B$ is hereditarily essential. Let $C$ be a unital intermediate $C^*$-algebra. Let $I$ be a non-trivial ideal in $C$, then $I \cap A \neq \{0\}$ implies $I \cap A = A$, hence $I$ contains the identity of $C$ and so $C$ is simple. On the other hand, suppose that $A \subset B$ is $C^*$-irreducible and let $C$ be an intermediate unital $C^*$-algebra. Then the only non-trivial ideal in $C$ is $C$ itself and so it intersects $A$ non-trivially. 
\end{proof}

\begin{prop}
    Let $A \subset B$ be an inclusion of unital $C^*$-algebras and suppose that there is a unital $C^*$-algebra $D \subset A$ such that every pseudo-conditional expectation from $B$ to $D$ is faithful. Then $A \subset B$ is hereditarily essential. In particular, if $A$ and $B$ are simple, the inclusion $A \subset B$ is $C^*$-irreducible.
\end{prop}
\begin{proof}
    By \cite{PiZa} Theorem 3.5 we need to show that every pseudo-expectation from $B$ to $A$ is faithful. Suppose this is not the case and let $\phi : B \rightarrow I(A)$ be a pseudo-expectation with non-trivial kernel. Let $\psi : I(A) \rightarrow I(D)$ be an extension of a pseudo-expectation from $A$ to $D$ (which exists by injectivity of $I(D)$). The composition $\psi \circ \phi$ is a non-faithful pseudo-expectation from $B$ to $D$, which is impossible. The result follows from Lemma \ref{lem1}.
\end{proof}
The same conclusion can be reached using \cite[Prop. 1.2.2]{Za}.

\begin{cor}
    The inclusion $\cO_2 \subset \cQ_2$ is $C^*$-irreducible.
\end{cor}
\begin{proof}
    The diagonal subalgebra $\cD_2 \subset \cO_2$ is Cartan in $\cQ_2$ \cite[Section 3]{ACR}. hence it admits a unique pseudo-expectation from $\cQ_2$ (by virtue of \cite{PiZa} Theorem 1.4), but the conditional expectation from $\cQ_2$ onto $\cD_2$ considered in \cite{LarsenLi} is faithful.
\end{proof}
Now consider the square of $C^*$-algebras
\[
\begin{array}{ccc}
\cO_2 & \subset & \cQ_2\\
\cup && \cup  \\
\cF_2 & \subset & \cB_2 \\
\end{array}
\]
where $\cF_2 = \cO_2^{\mathbb T}$ and $\cB_2 = \cQ_2^{\mathbb T}$ are the core UHF algebra and the Bunce-Deddens algebras of type $2^\infty$, respectively. Since $\cD_2 \subset \cF_2$ is Cartan in all of them, it readily follows, by the same argument, that all the horizontal, vertical and diagonal inclusions are $C^*$-irreducible as well.

\medskip
Several examples of irreducible, or more generally hereditarily essential inclusions of $C^*$-algebras are known. An interesting family comes from dynamical considerations. For example, arguing as in \cite{AmBa}, it is possible to see that for certain ``negatively curved'' groups there are abelian subgroups such that the inclusion of the corresponding reduced group $C^*$-algebras is hereditarily essential. Interestingly, this property never passes to the weak closures, denying the possibility to approach properties like solidity in this way (at least for the time being), the latter property being strictly related to hyperbolicity and the (AO)-property (\cite{ozas, AkOs, Ba, BaRa})


\begin{remark}
In \cite[Example 5.11]{Ro} it is observed that the inclusion $\cF_n
\subset \cO_n$ is $C^*$-irreducible, and further that any intermediate
$C^*$-algebra for this inclusion has the form $C^*(\cF_n, S_1^d) =
\cO_n^{{\mathbb Z}_d}$ for some integer $d \geq 2$, the fixed-point
algebra of $\cO_n$ under the order $d$ gauge automorphisms
$\lambda_{\omega 1}$, where $\omega \in {\mathbb T}$ is a primitive
$d$-root of $1$ and $\omega 1 \in \cU(\cO_n)$ (here we use the standard
terminology for endomorphisms of the Cuntz algebras, see e.g. \cite{ACR3} Section 2.2). Moreover, it is not
difficult to see that
$\cO_n^{{\mathbb Z}_d}$ is then isomorphic to $\cO_{n^d}$ (see e.g.
\cite[Prop. 7.2]{ACR3}). Similarly,
one can produce other intermediate $C^*$-algebras for the inclusion $\cF_n
\subset \cQ_n$ of the form $\cQ_n^{{\mathbb Z}_d}$, the fixed-point
algebra of $\cQ_n$ under the extended gauge automorphism
$\tilde\lambda_{\omega 1}$ mapping $U$ to $U$ and $S_n$ to $\omega S_n$,
and prove that
$\cQ_n^{{\mathbb Z}_d}$  is isomorphic to $\cQ_{n^d}$ (cf. \cite[Prop.
9.12]{ACR3} for the case $n=2$), where in general $\cQ_n \supset \cO_n$ is
the universal $C^*$-algebra generated by a unitary $V$ and an isometry
$S_n$ such that $V^n S_n = S_n V$ and $\sum_{k=0}^{n-1} V^k S_n S_n^*
V^{-k} = 1$. 
\end{remark}

Clearly, any intermediate $C^*$-subalgebra $\cE$ between $\cO_2$ and $\cQ_2$ in addition to being simple is also exact.
Moreover, it is possible to show with some additional work
 that $\cE$ must be 
purely infinite.  
One might wonder whether it has to be nuclear 
too. 
Anyway, one should stress that it is not clear at all if any such nontrivial $\cE$ actually exists.
A possible candidate for such an $\cE$ could be the $C^*$-algebra generated by the unitary normalizer
$N_{\cQ_2}(\cO_2)$, which is known to be strictly contained in $\cQ_2$ \cite[Theorem 7.4]{ACR2},
but it seems to be unknown whether $N_{\cQ_2}(\cO_2) = \cU(\cO_2)$.

\medskip
In \cite[Theorem 4.5]{ACR} it has been shown that if $\Lambda$ is a unital $*$-endomorphism of $\cQ_2$ such that $\Lambda|_{\cO_2} = {\rm id}_{\cO_2}$ then $\Lambda = {\rm id}_{\cQ_2}$. 
We can push it a bit further.
We say that a unital inclusion of $C^*$-algebras $A \subset B$ is {\it rigid} if ${\rm id}_B$ is the only ucp-map $\phi$ from $B$ into itself such that $\phi|_A = {\rm id}_A$, cf. \cite[Def. 4.3]{Su} (actually this definition can be traced back to the work of Hamana on injective envelopes). \\

Recall that the $C^*$-algebra $\cQ_2$ (and hence also $\cO_2$) is represented canonically (and faithfully) on $\ell^2(\mathbb{Z})$ (\cite{LarsenLi} Section 2).

\begin{prop}
\label{prop2.6}
Let $\cQ_2$ and $\cO_2$ be represented canonically on $H=\ell^2(\mathbb{Z})$. Every ucp map $\cQ_2 \rightarrow \mathbb{B}(H)$ which is the identity on $\cO_2$ is the identity on $\cQ_2$. In particular, the inclusion $\cO_2 \subset \cQ_2$ is rigid.
\end{prop}
\begin{proof}
Let $\phi : \cQ_2 \rightarrow \mathbb{B}(H)$ be a ucp map such that its restriction to $\cO_2$ is the identity. For every $k \in \mathbb{N}$ we define the projections
\begin{equation}
    p_k:=S_1^k S_2 S_2^* (S_1^k)^*, \ k \in {\mathbb N}
\end{equation}
Since $\cO_2$ is in the multiplicative domain of $\phi$ we have, using the relation $U S_1^k S_2 = S_2^k S_1$ , that for every $k$
\begin{equation}
\phi (US_1^k S_2) 
= \phi(U)S_1^k S_2
= \phi (S_2^k S_1)= S_2^k S_1= US_1^k S_2,
\end{equation}
from which we obtain
\begin{equation}
    \phi (Up_k)=\phi (U) p_k = U p_k.
\end{equation}
Since $\sum_{k=1}^\infty p_k $ converges weakly to $1$ in the canonical representation of $\cQ_2$ on $\ell^2 (\mathbb{Z})$, it readily follows that $\phi(U)=U$. 
It is now routine to verify that $\phi$ is the identity map, as desired. 
\end{proof}

\begin{cor}
\label{corinj}
The injective envelope of $\cO_2$ is $*$-isomorphic to the injective envelope of $\cQ_2$.
\end{cor}
\begin{proof}
    Realize both the injective envelope of $\cO_2$ and the injective envelope of $\cQ_2$ inside $\mathbb{B}(H)$, where $H =\ell^2(\mathbb{Z})$ (for example using \cite{ham} Theorem 3.4). By injectivity, the inclusion  $\cO_2 \subset I(\cO_2)$ extends to a ucp map $\psi : I(\cQ_2) \rightarrow I(\cO_2)$ and the inclusion $\cO_2 \subset I(\cQ_2)$ extends to a ucp map $\phi: I(\cO_2) \rightarrow I(\cQ_2)$. The restriction of $\psi$ to $\cO_2$ is the identity, hence by Proposition \ref{prop2.6} it is the identity on $\cQ_2$. In particular $\cO_2 \subset \cQ_2 \subset I(\cO_2)$. Now both $\phi \circ \psi$ and $\psi \circ \phi$ restrict to the identity on $\cQ_2$ (hence also on $\cO_2$); by the universal property of the injective envelope (see \cite{ham} Definition 2.2 and Theorem 4.1) $\phi \circ \psi =\id_{I(\cO_2)}$ and $\psi \circ \phi =\id_{I(\cQ_2)}$. It follows that both $\phi$ and $\psi$ are injective, hence they are complete order isomorphisms and thus $*$-isomorphisms (\cite{blackadar} Theorem II.6.9.17).
\end{proof}

The authors 
are not aware of other examples of inclusions of $C^*$-algebras for which the conclusion of Corollary \ref{corinj} holds. 
Inclusions of $C^*$-algebras sharing the same injective envelope are abundant in the equivariant setting and actually this fact is relevant in the 
discussion of a conjecture by Ozawa (\cite{Oz2,KaKe,BaRa2})
\medskip

Summing up, we have shown the following result.
\begin{theorem}
The natural inclusion $\cO_2 \subset \cQ_2$ satisfies the following properties:
\begin{itemize}
\item it is $C^*$-irreducible;
\item it is rigid. 
\end{itemize}
Moreover, the injective envelopes of $\cO_2$ and $\cQ_2$ are isomorphic.
\end{theorem}

We have discussed the case $n=2$ since it is easier to grasp relevant results from the existing literature. 
It is likely that similar results 
carry over to the case of $\cO_n \subset \cQ_n$ for all $n >2$.

\section{Acknowledgments}
Both authors thank the anonymous referee for the precious comments which led to an improved exposition. The first named author ackowledges the support of INdAM-GNAMPA, IM PAN, the grant SOE-Young Researchers 2024 "Generalized Akemann-Ostrand Property: Analytical, Dynamical And Rigidity Properties (GAOPADRP)", CUP: E83C24002550001 and the MIUR Excellence Department Project 2023--2027 MatMod@Tov awarded to the Department of Mathematics, University of Rome Tor Vergata.

\bigskip
{\parindent=0pt Addresses of the authors:

\medskip


Jacopo Bassi, \textsc{Department of Mathematics, University of Tor Vergata, Via della Ricerca Scientifica 1, 00133 Roma, Italy }\par\nopagebreak
E-mail: \texttt{bassi@axp.mat.uniroma2.it}
\par}

\bigskip \noindent
Roberto Conti, \textsc{Dipartimento SBAI,
Sapienza Universit\`a di Roma, 
Via A. Scarpa 16,
I-00161 Roma, Italy.}
\\ E-mail: \texttt{roberto.conti@sbai.uniroma1.it}
\par

\end{document}